# Adaptive geometric multigrid for the mixed finite cell formulation of Stokes and Navier-Stokes equations [*]

S. Saberi [†]    G. Meschke [‡]    A. Vogel [†]

**Abstract**

Unfitted finite element methods have emerged as a popular alternative to classical finite element methods for the solution of partial differential equations and allow modeling arbitrary geometries without the need for a boundary-conforming mesh. On the other hand, the efficient solution of the resultant system is a challenging task because of the numerical ill-conditioning that typically entails from the formulation of such methods. We use an adaptive geometric multigrid solver for the solution of the mixed finite cell formulation of saddle-point problems and investigate its convergence in the context of the Stokes and Navier-Stokes equations. We present two smoothers for the treatment of cutcells in the finite cell method and analyze their effectiveness for the model problems using a numerical benchmark. Results indicate that the presented multigrid method is capable of solving the model problems independently of the problem size and is robust with respect to the depth of the grid hierarchy.

## 1  Introduction

The finite element method is a powerful tool for the numerical approximation of partial differential equations and has been successfully applied to a wide range of problems in the past decades. The generation of an appropriate tessellation of the computational domain often is a daunting task in boundary-conforming finite element methods and has motivated the development of a number of techniques, such as the extended finite element method (XFEM) [4], cutFEM [9] and the finite cell method [34, 14, 45] under the umbrella of unfitted finite element methods. The mutual goal of such methods is to separate the physical domain from the computational domain so that an arbitrary geometry can be resolved by a trivially constructed computational mesh, circumventing the need for the generation of a boundary-conforming tessellation of the domain. The finite cell method is a fictitious domain method and employs an integration technique such as Gaussian over-integration and uniform or adaptive subcell integration for the resolution of the physical domain and a method for the weak imposition of essential boundary conditions such as Lagrange multipliers [19, 20, 22, 10], the penalty method [2, 51, 8] and the Nitsche's method [33, 24, 13, 18, 11]. The finite cell method has been applied to different classes of problems, including structural mechanics [34, 14] and fluid dynamics [25], see [45] for a review.

[*]Financial support was provided by the German Research Foundation (*Deutsche Forschungsgemeinschaft, DFG*) in the framework of the collaborative research center SFB 837 *Interaction Modeling in Mechanized Tunneling*.
[†]High Performance Computing in the Engineering Sciences, Ruhr University Bochum, Universitätsstr. 150, 44801 Bochum, Germany
[‡]Institute for Structural Mechanics, Ruhr University Bochum, Universitätsstr. 150, 44801 Bochum, Germany

A wide variety of problems in fluid dynamics, structural mechanics, finance, etc. are mathematically described as saddle-point problems. Several methods, including Uzawa methods [17, 7], multigrid methods [49, 23, 50] and Krylov subspace solvers [6, 40, 47, 41] are typically employed for the iterative solution of saddle-point systems, see also [16, 37, 32]. Effective preconditioning is essential for Krylov subspace solvers to achieve scalable convergence. Multigrid methods are often used as preconditioners in Krylov subspace solvers. Other preconditioners include Schur complement methods, such as pressure correction schemes [35], approximate commutator schemes [46, 28, 15], and domain decomposition methods, such as overlapping Schwarz methods [30, 31, 36].

Geometric multigrid methods are among the most efficient algorithms for the solution of saddle-point problems [16, 44]. A multigrid method based on incomplete LU factorization was proposed in [50]. The performance of some pressure correction methods as smoothers was studied in [21]. A pressure correction smoother, known as Braess-Sarazin smoother, was proposed in [5]. The Vanka smoother, which is based on the solution of local saddle-point problems, was proposed in [48].

The iterative solution of saddle-point problems is a challenging task which is further complicated by the fact that the weak formulation of the finite cell method typically leads to matrices that are ill-conditioned. The conditioning of the system matrix depends on how the physical domain intersects the computational mesh and deteriorates in the presence of small cut fractions. A conjugate gradient method for the finite cell formulation of the elasticity problem was studied in [27] with an AMG preconditioner and in [26] with an additive Schwarz preconditioner. A preconditioner for the conjugate gradient and GMRES methods was studied for the solution of a variety of finite cell problems in [38]. Geometric multigrid with Schwarz-type smoothers was recently studied for the finite cell formulation of elliptic problems in [39, 42].

We use the mixed finite cell formulation of the incompressible Stokes and Navier-Stokes equations with Nitsche's boundary conditions. We employ space trees [12] for the spatial discretization of the computational domain and use adaptive integration in addition to adaptive mesh refinement. The contributions of this work can be summarized as follows:

- We formulate an adaptive geometric multigrid solver for the solution of the mixed finite cell formulation of saddle-point problems in the context of the Stokes and Navier-Stokes equations and present two smoothers for the treatment of the resultant systems.

- We investigate the convergence of the geometric multigrid method in different scenarios both as a solver and as a preconditioner in Krylov subspace solvers and evaluate its performance with respect to the choice of smoothers using numerical benchmarks for the model problems.

The remainder of this paper is organized as follows. The mixed finite cell formulation of the incompressible Stokes and Navier-Stokes equations is derived in Section 2. The necessary components of the geometric multigrid solver, including the generation of the grid hierarchy and smoothers are discussed in Section 3. The performance of the developed multigrid method is analyzed in Section 4. Finally, some concluding remarks are given in Section 5.

## 2 Mixed finite cell formulation of the model problems

Given a physical domain $\Omega$ with boundary $\partial\Omega$, we consider a regular embedding domain $\Omega_e$ such that $\Omega \subset \Omega_e$, see Figure 1(a). At the discrete level, $\Omega_{e,h}$ is composed of a set of $n_c$ cells $\mathcal{S}_c := \{c_i \mid c_i \in \Omega_{e,h}\}$ with a set of $n_{cc}$ cut cells $\mathcal{S}_{cc} := \{c_i \in \Omega_{e,h} \mid c_i \cap \partial\Omega \neq 0\}$. We start by deriving the mixed finite cell formulation of the incompressible Stokes equations. The strong form of the Stokes equations can be written as

$$\begin{aligned}-\eta\nabla^2\boldsymbol{u} + \nabla p &= \boldsymbol{f} \quad \text{in } \Omega, \\ \nabla \cdot \boldsymbol{u} &= 0 \quad \text{in } \Omega,\end{aligned} \quad (1)$$

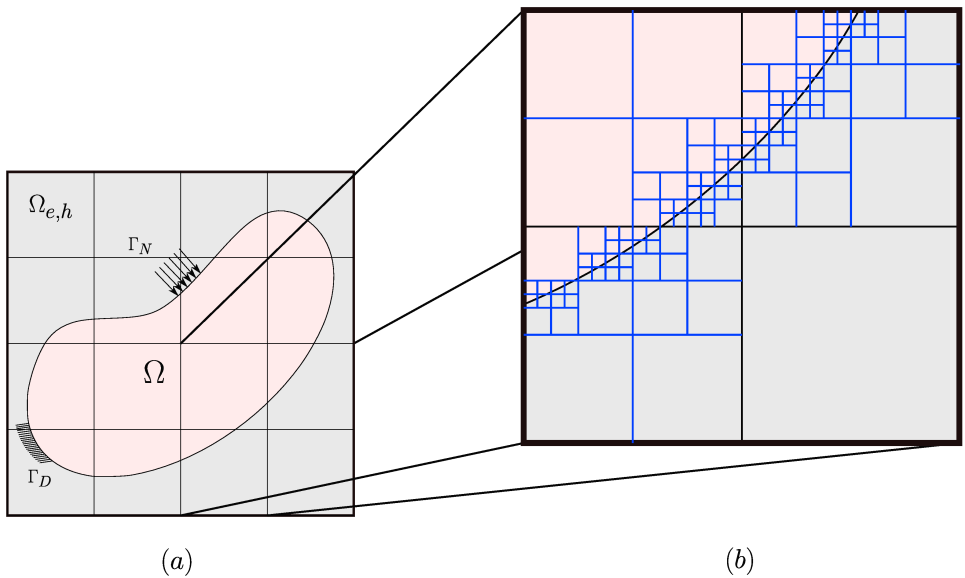

(a)             (b)

Figure 1: (a) The physical domain $\Omega$, the embedding domain $\Omega_{e,h}$ with a uniform discretization and Dirichlet $\Gamma_D$ and (inhomogeneous) Neumann $\Gamma_N$ boundary conditions and (b) adaptive integration of a part of the domain. The blue lines designate subcells, adaptively generated towards the boundary of the physical domain, and only appear during the integration stage, i.e., do not exist in the discrete system

where $\boldsymbol{u}$ is the fluid velocity vector, $p$ is the fluid pressure, $\boldsymbol{f}$ is the body force exerted on the fluid and $\eta$ is the kinematic viscosity.

We denote by $(\cdot,\cdot)_\Omega$ and $(\cdot,\cdot)_{\partial\Omega}$ the $L^2$ scalar product on $\Omega$ and $\partial\Omega$, respectively. The weak form of the Stokes equations is obtained as follows by multiplying Equation 1 by test functions from the left, integrating over the computational domain and transferring the derivative from the pressure trial function to the velocity test functions using integration by parts

$$(\eta\nabla\boldsymbol{v},\nabla\boldsymbol{u})_\Omega - (\boldsymbol{v},\boldsymbol{n}\cdot\eta\nabla\boldsymbol{u})_{\partial\Omega} - (\nabla\cdot\boldsymbol{v},p)_\Omega + (\boldsymbol{v},\boldsymbol{n}p)_{\partial\Omega} - (q,\nabla\cdot\boldsymbol{u})_\Omega = \\ (\boldsymbol{v},\boldsymbol{f})_\Omega, \quad (2)$$

where $(\boldsymbol{v},q)$ are the infinite-dimensional vector-valued velocity and scalar pressure test functions, respectively. The Stokes boundary value problem is formed by requiring the solution to Equation 1 to satisfy the following boundary conditions

$$\boldsymbol{u} = \boldsymbol{w} \text{ on } \Gamma_D \subset \partial\Omega,$$
$$\eta\frac{\partial\boldsymbol{u}}{\partial\boldsymbol{n}} - \boldsymbol{n}p = \boldsymbol{h} \text{ on } \Gamma_N := \partial\Omega \setminus \Gamma_D, \quad (3)$$

where $\partial\Omega$ is the boundary of the computational domain, $\Gamma_D$ and $\Gamma_N$ denote the Dirichlet and Neumann parts of the boundary, respectively, such that $\partial\Omega = \Gamma_D \cup \Gamma_N$ and $\Gamma_D \cap \Gamma_N = 0$ and $\boldsymbol{n}$ is the outward unit normal vector to the boundary.

The spaces of test and trial functions in the weak form of the Stokes equations in boundary-conforming finite element methods are chosen as

$$\begin{aligned}
\boldsymbol{V}_w &:= \{\boldsymbol{u} \in H^1(\Omega)^d \mid \boldsymbol{u} = \boldsymbol{w} \text{ on } \Gamma_D\}, \\
\boldsymbol{V}_0 &:= \{\boldsymbol{v} \in H^1(\Omega)^d \mid \boldsymbol{v} = \boldsymbol{0} \text{ on } \Gamma_D\}, \\
Q &:= \{q \in L^2(\Omega)^d\},
\end{aligned} \quad (4)$$

where $d$ is the dimension of the domain. The weak form in Equation 2 can be written as follows by combining the boundary terms and separating the boundary into its Dirichlet and Neumann parts:

Find $(\boldsymbol{u}, p) \in (\boldsymbol{V}_w, Q)$ such that

$$(\eta \nabla \boldsymbol{v}, \nabla \boldsymbol{u})_\Omega - (\nabla \cdot \boldsymbol{v}, p)_\Omega - (q, \nabla \cdot \boldsymbol{u})_\Omega = \\ (\boldsymbol{v}, \boldsymbol{f})_\Omega + (\boldsymbol{v}, \boldsymbol{h})_{\Gamma_N}, \quad \forall (\boldsymbol{v}, q) \in (\boldsymbol{V}_0, Q) \tag{5}$$

The imposition of essential boundary conditions in the form above is referred to as the strong imposition of boundary conditions. In unfitted finite element methods, however, the strong imposition of essential boundary conditions is in general not possible as the physical domain does not necessarily coincide with the computational domain. Therefore, essential boundary conditions are imposed weakly using methods such as Lagrange multipliers [19, 20, 22, 10], the penalty method [2, 51, 8] and Nitsche's method [33, 24, 13, 18, 11]. The trial and test functions are now chosen from the same spaces $(\boldsymbol{V}, Q)$, given by

$$\boldsymbol{V} := \{\boldsymbol{u} \in H^1(\Omega)^d\}, \\ Q := \{q \in L^2(\Omega)^d\}. \tag{6}$$

Therefore, the boundary terms in Equation 2 do not vanish over the Dirichlet boundary. We employ Nitsche's method for the imposition of essential boundary conditions in this work. The weak form of the Stokes equations with Nitsche's boundary terms is written as

$$(\eta \nabla \boldsymbol{v}, \nabla \boldsymbol{u})_\Omega - (\boldsymbol{v}, \boldsymbol{n} \cdot \eta \nabla \boldsymbol{u})_{\Gamma_D} - (\boldsymbol{n} \cdot \eta \nabla \boldsymbol{v}, \boldsymbol{u} - \boldsymbol{w})_{\Gamma_D} + (\boldsymbol{v}, \lambda(\boldsymbol{u} - \boldsymbol{w}))_{\Gamma_D} \\ -(\nabla \cdot \boldsymbol{v}, p)_\Omega + (\boldsymbol{v}, \boldsymbol{n}p)_{\Gamma_D} + (\boldsymbol{n}q, (\boldsymbol{u} - \boldsymbol{w}))_{\Gamma_D} - (q, \nabla \cdot \boldsymbol{u})_\Omega = \\ (\boldsymbol{v}, \boldsymbol{f})_\Omega + (\boldsymbol{v}, \boldsymbol{h})_{\Gamma_N}. \tag{7}$$

*Remark* 1. The third and seventh terms on the l.h.s of Equation 7 are the velocity and pressure symmetric consistency terms, respectively. These terms ensure that the modified weak form retains the symmetry of the original weak form. The fourth term on the l.h.s of Equation 7 is a stabilization term with a positive stabilization parameter $\lambda$ that ensures the stability of the solution for a large-enough $\lambda$ and that the solution satisfies the essential boundary conditions. These terms also retain the variational consistency in Nitsche's method (see Remark 2).

The boundary terms in Equation 7 can be separated to form the following weak problem:

Find $(\boldsymbol{u}, p) \in (\boldsymbol{V}, Q)$ such that

$$(\eta \nabla \boldsymbol{v}, \nabla \boldsymbol{u})_\Omega - (\boldsymbol{v}, \boldsymbol{n} \cdot \eta \nabla \boldsymbol{u})_{\Gamma_D} - (\boldsymbol{n} \cdot \eta \nabla \boldsymbol{v}, \boldsymbol{u})_{\Gamma_D} + (\boldsymbol{v}, \lambda \boldsymbol{u})_{\Gamma_D} \\ -(\nabla \cdot \boldsymbol{v}, p)_\Omega + (\boldsymbol{v}, \boldsymbol{n}p)_{\Gamma_D} + (\boldsymbol{n}q, \boldsymbol{u})_{\Gamma_D} - (q, \nabla \cdot \boldsymbol{u})_\Omega = \\ (\boldsymbol{v}, \boldsymbol{f})_\Omega + (\boldsymbol{v}, \boldsymbol{h})_{\Gamma_N} - (\boldsymbol{n} \cdot \eta \nabla \boldsymbol{v}, \boldsymbol{w})_{\Gamma_D} + (\boldsymbol{v}, \lambda \boldsymbol{w})_{\Gamma_D} + (\boldsymbol{n}q, \boldsymbol{w})_{\Gamma_D}, \\ \forall (\boldsymbol{v}, q) \in (\boldsymbol{V}, Q). \tag{8}$$

*Remark* 2. The weak form in Equation 8 is consistent in the sense that it can be shown that the solution to the original problem in Equation 1 with the boundary conditions in Equation 3 also satisfies the weak form in Equation 8.

*Remark* 3. Compared to the Lagrange multiplier method and the penalty method, the Nitsche's method has the advantage of keeping the problem size the same as the original problem, retaining the symmetry of the original weak form and being variationally consistent.

We now derive the finite cell weak form of the Stokes equations with Nitsche's boundary conditions by extending the weak form in Equation 8 to an embedding domain as follows:

Find $(\boldsymbol{u}, p) \in (\boldsymbol{V}_e, Q_e)$ such that

$$(\eta \nabla \boldsymbol{v}, \alpha \nabla \boldsymbol{u})_{\Omega_e} - (\boldsymbol{v}, \boldsymbol{n} \cdot \eta \nabla \boldsymbol{u})_{\Gamma_D} - (\boldsymbol{n} \cdot \eta \nabla \boldsymbol{v}, \boldsymbol{u})_{\Gamma_D} + (\boldsymbol{v}, \lambda \boldsymbol{u})_{\Gamma_D} \\ -(\nabla \cdot \boldsymbol{v}, \alpha p)_{\Omega_e} + (\boldsymbol{v}, \boldsymbol{n}p)_{\Gamma_D} + (\boldsymbol{n}q, \boldsymbol{u})_{\Gamma_D} - (q, \alpha \nabla \cdot \boldsymbol{u})_{\Omega_e} = \\ (\boldsymbol{v}, \alpha \boldsymbol{f})_{\Omega_e} + (\boldsymbol{v}, \boldsymbol{h})_{\Gamma_N} - (\boldsymbol{n} \cdot \eta \nabla \boldsymbol{v}, \boldsymbol{w})_{\Gamma_D} + (\boldsymbol{v}, \lambda \boldsymbol{w})_{\Gamma_D} + (\boldsymbol{n}q, \boldsymbol{w})_{\Gamma_D}, \\ \forall (\boldsymbol{v}, q) \in (\boldsymbol{V}_e, Q_e), \tag{9}$$

where $\Omega_e$ is the embedding domain, $(\boldsymbol{V}_e, Q_e)$ are the spaces in Equation 6 defined over $\Omega_e$ and $\alpha$ is a penalization parameter chosen as

$$\begin{cases} \alpha = 1 & \text{in } \Omega, \\ \alpha = 0 & \text{in } \Omega_e \setminus \Omega. \end{cases} \tag{10}$$

*Remark* 4. It can be seen, that with this choice of the penalization parameter, the weak form in Equation 9 recovers the weak form in Equation 8.

Finally, introducing discrete spaces $(\boldsymbol{V}_{e,h}, Q_{e,h}) \subset (\boldsymbol{V}_e, Q_e)$, the discrete bilinear and linear forms of the mixed finite cell formulation of the Stokes equations can be written as

$$\begin{aligned} a(\boldsymbol{v}_h, \boldsymbol{u}_h) + b(\boldsymbol{v}_h, p_h) &= f(\boldsymbol{v}_h) \\ b(q_h, \boldsymbol{u}_h) &= g(q_h), \end{aligned} \tag{11}$$

with

$$\begin{aligned} a(\boldsymbol{v}_h, \boldsymbol{u}_h) &= (\eta \nabla \boldsymbol{v}_h, \alpha \nabla \boldsymbol{u}_h)_{\Omega_{e,h}} - (\boldsymbol{v}_h, \boldsymbol{n} \cdot \eta \nabla \boldsymbol{u}_h)_{\Gamma_{D,h}} - (\boldsymbol{n} \cdot \eta \nabla \boldsymbol{v}_h, \boldsymbol{u}_h)_{\Gamma_{D,h}} \\ &\quad + (\boldsymbol{v}_h, \lambda \boldsymbol{u}_h)_{\Gamma_{D,h}} \\ b(\boldsymbol{v}_h, p_h) &= -(\nabla \cdot \boldsymbol{v}_h, \alpha p_h)_{\Omega_{e,h}} + (\boldsymbol{v}_h, \boldsymbol{n} p_h)_{\Gamma_{D,h}} \\ f(\boldsymbol{v}_h) &= (\boldsymbol{v}_h, \alpha \boldsymbol{f})_{\Omega_{e,h}} + (\boldsymbol{v}_h, \boldsymbol{h})_{\Gamma_{N,h}} - (\boldsymbol{n} \cdot \eta \nabla \boldsymbol{v}_h, \boldsymbol{w})_{\Gamma_{D,h}} + (\boldsymbol{v}_h, \lambda \boldsymbol{w})_{\Gamma_{D,h}} \\ g(q_h) &= (\boldsymbol{n} q_h, \boldsymbol{w})_{\Gamma_D}, \end{aligned} \tag{12}$$

where $\Omega_{e,h}$ and $\Gamma_{\cdot,h}$ are appropriate discretizations of the embedding domain and the physical boundary, respectively.

*Remark* 5. The penalization parameter $\alpha$ in the discrete bilinear and linear forms in Equation 11 is taken as a small value ($\alpha \ll 1$) instead of absolute zero outside of the physical domain in order to mitigate the numerical ill-conditioning of the bilinear form.

The discrete form in Equation 11 leads to a linear system of equations of the form

$$\begin{bmatrix} \boldsymbol{A} & \boldsymbol{B} \\ \boldsymbol{B}^T & \boldsymbol{C} \end{bmatrix} \begin{bmatrix} \boldsymbol{u} \\ \boldsymbol{p} \end{bmatrix} = \begin{bmatrix} \boldsymbol{f} \\ \boldsymbol{g} \end{bmatrix}, \tag{13}$$

where the $\boldsymbol{A}$, $\boldsymbol{B}$, $\boldsymbol{f}$ and $\boldsymbol{g}$ are defined according to Equation 12 and $\boldsymbol{C}$ is either zero in the case of conforming spaces that satisfy the inf-sup condition or the stabilization block otherwise. We employ a Q1-Q1 discretization, where the stabilization term is defined as

$$\boldsymbol{C} := -\beta \sum_{\Omega_{c_i} \in \Omega_{e,h}} h_{\Omega_{c_i}}^2 (\nabla \boldsymbol{v}_h, \nabla \boldsymbol{u}_h)_{\Omega_{c_i}}, \tag{14}$$

where $\beta$ is a sufficiently large positive constant and $\Omega_{c_i}$ is the domain of $c_i$.

We now focus on the mixed finite cell formulation of the incompressible Navier-Stokes equations. Navier-Stokes equations are nonlinear equations that can describe flow regimes outside of the scope of Stokes equations, e.g., at high Reynold numbers. The strong form of the Navier-Stokes equations is obtained by adding a nonlinear convection term to the Stokes equations (see Equation 1) and is given by

$$\begin{aligned} -\eta \nabla^2 \boldsymbol{u} + \boldsymbol{u} \cdot \nabla \boldsymbol{u} + \nabla p &= \boldsymbol{f} \quad \text{in } \Omega, \\ \nabla \cdot \boldsymbol{u} &= 0 \quad \text{in } \Omega. \end{aligned} \tag{15}$$

Equation 15 along with the boundary conditions in Equation 3 form the Navier-Stokes boundary value problem. The mixed finite cell weak form of the Navier-Stokes equations is obtained following the same procedure described above for the Stokes equations:

Find $(\boldsymbol{u}, p) \in (\boldsymbol{V}_e, Q_e)$ such that

$$\begin{aligned} F(\boldsymbol{u}, p) &:= (\eta \nabla \boldsymbol{v}, \alpha \nabla \boldsymbol{u})_{\Omega_e} - (\boldsymbol{v}, \boldsymbol{n} \cdot \eta \nabla \boldsymbol{u})_{\Gamma_D} - (\boldsymbol{n} \cdot \eta \nabla \boldsymbol{v}, \boldsymbol{u})_{\Gamma_D} + (\boldsymbol{v}, \lambda \boldsymbol{u})_{\Gamma_D} \\ &\quad + (\boldsymbol{v}, \boldsymbol{u} \cdot \nabla \boldsymbol{u})_{\Omega_e} - (\nabla \cdot \boldsymbol{v}, \alpha p)_{\Omega_e} + (\boldsymbol{v}, \boldsymbol{n} p)_{\Gamma_D} + (\boldsymbol{n} q, \boldsymbol{u})_{\Gamma_D} - (q, \alpha \nabla \cdot \boldsymbol{u})_{\Omega_e} \\ &\quad - (\boldsymbol{v}, \alpha \boldsymbol{f})_{\Omega_e} - (\boldsymbol{v}, \boldsymbol{h})_{\Gamma_N} + (\boldsymbol{n} \cdot \eta \nabla \boldsymbol{v}, \boldsymbol{w})_{\Gamma_D} - (\boldsymbol{v}, \lambda \boldsymbol{w})_{\Gamma_D} - (\boldsymbol{n} q, \boldsymbol{w})_{\Gamma_D} = 0, \\ &\quad \forall (\boldsymbol{v}, q) \in (\boldsymbol{V}_e, Q_e). \end{aligned} \tag{16}$$

The weak form of the Navier-Stokes equations is solved using nonlinear iterative methods, where the iterations comprise a series of linearized problems. We discuss two linearization methods, namely the Newton's method and the Picard's method. Let $G(\boldsymbol{u}, p) : (\boldsymbol{V}_e, Q_e) \to (\boldsymbol{V}_e, Q_e)$ be a function such that $G(\boldsymbol{u}^*, p^*) = (\boldsymbol{u}^*, p^*)$ if $F(\boldsymbol{u}^*, p^*) = 0$. The iterations $(\boldsymbol{u}^{k+1}, p^{k+1}) = G(\boldsymbol{u}^k, p^k)$ can then be shown to converge to $(\boldsymbol{u}^*, p^*)$ if $G$ is a contraction [29].

The Newton iteration can then be formulated as $G_{\text{Newton}} := (\boldsymbol{u}^k, p^k) + \frac{R(\boldsymbol{u}^k, p^k)}{F'(\boldsymbol{u}^k, p^k)}$, where $R$ is the residual and $F'$ is the Fréchet derivative of the weak form in Equation 16. The Newton iteration is commonly formulated as a function of $(\delta \boldsymbol{u}^k, \delta p^k) := (\boldsymbol{u}^{k+1}, p^{k+1}) - (\boldsymbol{u}^k, p^k)$ as

$$F'(\boldsymbol{u}^k, p^k)\, (\delta \boldsymbol{u}^k, \delta p^k) = R(\boldsymbol{u}^k, p^k), \tag{17}$$

with

$$\begin{aligned} F'(\boldsymbol{u}^k, p^k)\, (\delta \boldsymbol{u}^k, \delta p^k) := &(\eta \nabla \boldsymbol{v}, \alpha \nabla \delta \boldsymbol{u})_{\Omega_e} - (\boldsymbol{v}, \boldsymbol{n} \cdot \eta \nabla \delta \boldsymbol{u})_{\Gamma_D} - (\boldsymbol{n} \cdot \eta \nabla \boldsymbol{v}, \delta \boldsymbol{u})_{\Gamma_D} \\ &+ (\boldsymbol{v}, \lambda \delta \boldsymbol{u})_{\Gamma_D} + (\boldsymbol{v}, \delta \boldsymbol{u} \cdot \nabla \boldsymbol{u}^k)_{\Omega_e} + (\boldsymbol{v}, \boldsymbol{u}^k \cdot \nabla \delta \boldsymbol{u})_{\Omega_e} - (\nabla \cdot \boldsymbol{v}, \alpha \delta p)_{\Omega_e} + (\boldsymbol{v}, \boldsymbol{n} \delta p)_{\Gamma_D} \\ &+ (\boldsymbol{n} q, \delta \boldsymbol{u})_{\Gamma_D} - (q, \alpha \nabla \cdot \delta \boldsymbol{u})_{\Omega_e}, \end{aligned} \tag{18}$$

and

$$\begin{aligned} R(\boldsymbol{u}^k, p^k) := &(\boldsymbol{v}, \alpha \boldsymbol{f})_{\Omega_e} + (\boldsymbol{v}, \boldsymbol{h})_{\Gamma_N} - (\boldsymbol{n} \cdot \eta \nabla \boldsymbol{v}, \boldsymbol{w})_{\Gamma_D} + (\boldsymbol{v}, \lambda \boldsymbol{w})_{\Gamma_D} + (\boldsymbol{n} q, \boldsymbol{w})_{\Gamma_D} \\ &- (\eta \nabla \boldsymbol{v}, \alpha \nabla \boldsymbol{u}^k)_{\Omega_e} + (\boldsymbol{v}, \boldsymbol{n} \cdot \eta \nabla \boldsymbol{u}^k)_{\Gamma_D} + (\boldsymbol{n} \cdot \eta \nabla \boldsymbol{v}, \boldsymbol{u}^k)_{\Gamma_D} - (\boldsymbol{v}, \lambda \boldsymbol{u}^k)_{\Gamma_D} \\ &- (\boldsymbol{v}, \boldsymbol{u}^k \cdot \nabla \boldsymbol{u}^k)_{\Omega_e} + (\nabla \cdot \boldsymbol{v}, \alpha p^k)_{\Omega_e} - (\boldsymbol{v}, \boldsymbol{n} p^k)_{\Gamma_D} - (\boldsymbol{n} q, \boldsymbol{u}^k)_{\Gamma_D} + (q, \alpha \nabla \cdot \boldsymbol{u}^k)_{\Omega_e} \end{aligned} \tag{19}$$

The Picard iteration is defined as $G_{\text{Picard}} := \frac{1}{\mathcal{P}} \mathcal{K}$ and is commonly formulated as a direct iteration for $(\boldsymbol{u}^{k+1}, p^{k+1})$ as

$$\mathcal{P}\, (\boldsymbol{u}^{k+1}, p^{k+1}) = \mathcal{K} \tag{20}$$

with

$$\begin{aligned} \mathcal{P}\, (\boldsymbol{u}^{k+1}, p^{k+1}) := &(\eta \nabla \boldsymbol{v}, \alpha \nabla \boldsymbol{u}^{k+1})_{\Omega_e} - (\boldsymbol{v}, \boldsymbol{n} \cdot \eta \nabla \boldsymbol{u}^{k+1})_{\Gamma_D} - (\boldsymbol{n} \cdot \eta \nabla \boldsymbol{v}, \boldsymbol{u}^{k+1})_{\Gamma_D} \\ &+ (\boldsymbol{v}, \lambda \boldsymbol{u}^{k+1})_{\Gamma_D} + (\boldsymbol{v}, \boldsymbol{u}^k \cdot \nabla \boldsymbol{u}^{k+1})_{\Omega_e} - (\nabla \cdot \boldsymbol{v}, \alpha p^{k+1})_{\Omega_e} + (\boldsymbol{v}, \boldsymbol{n} p^{k+1})_{\Gamma_D} \\ &+ (\boldsymbol{n} q, \boldsymbol{u}^{k+1})_{\Gamma_D} - (q, \alpha \nabla \cdot \boldsymbol{u}^{k+1})_{\Omega_e}, \end{aligned} \tag{21}$$

and

$$\mathcal{K} := -(\boldsymbol{v}, \alpha \boldsymbol{f})_{\Omega_e} - (\boldsymbol{v}, \boldsymbol{h})_{\Gamma_N} + (\boldsymbol{n} \cdot \eta \nabla \boldsymbol{v}, \boldsymbol{w})_{\Gamma_D} - (\boldsymbol{v}, \lambda \boldsymbol{w})_{\Gamma_D} - (\boldsymbol{n} q, \boldsymbol{w})_{\Gamma_D}. \tag{22}$$

*Remark* 6. The Picard's method linearizes Equation 16 by fixing the convection coefficient at the current velocity and can also be derived from the Newton's method by dropping the term $(\boldsymbol{v}, \delta \boldsymbol{u} \cdot \nabla \boldsymbol{u}^k)_{\Omega_e}$ in Equation 18.

*Remark* 7. The Newton's method can be shown to achieve quadratic convergence when the current iteration becomes sufficiently close to the solution, while the Picard's method achieves linear convergence in the general case. On the other hand, the sphere of convergence of the Picard's method is significantly larger compared to the Newton's method, i.e., the Picard's method is more likely to converge to the solution for an initial guess that is far from the solution, see also [29].

# 3  Adaptive geometric multigrid

Multigrid methods are typically employed in one of two fashions for the solution of saddle-point problems: in the first approach, multigrid is used to treat the entire system, and the coupling between velocity and pressure is therefore preserved. The second approach is to use a variant of Schur complement schemes and use multigrid to solve the resultant systems associated with the velocity block and the Schur complement system. In general, the coupling between velocity and pressure is at least partially lost in the latter approach. In this work, we focus on the former approach and develop a monolithic adaptive geometric multigrid method for the mixed finite cell formulation of saddle-point problems in the context of Stokes and Navier-Stokes systems.

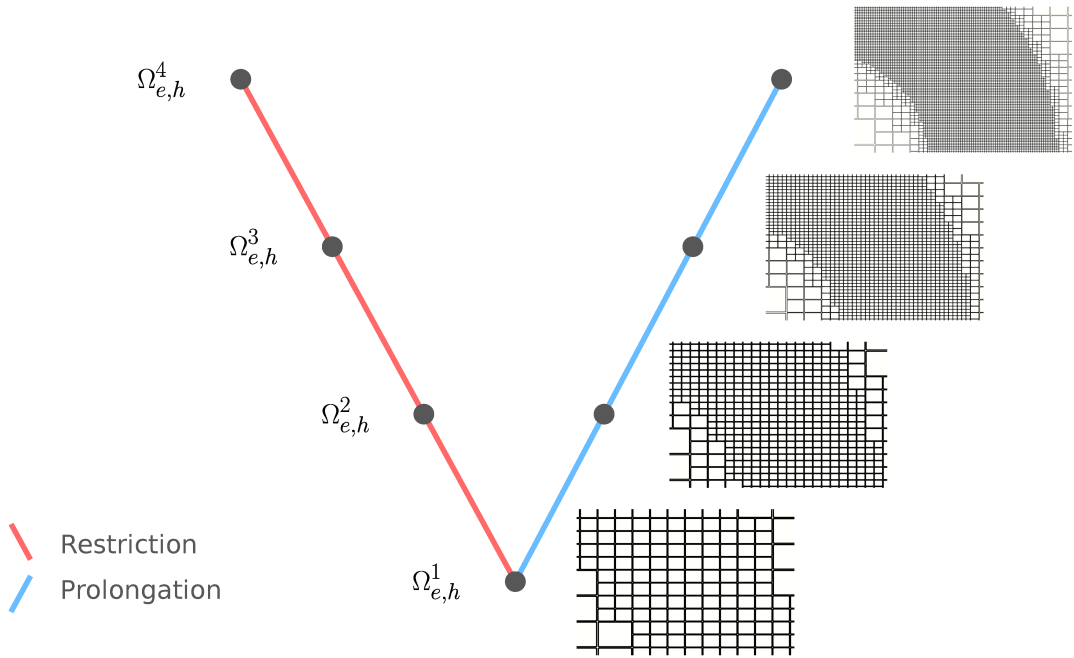

Figure 2: A four-level multigrid V cycle along with the adaptive refinement of a sample domain around an arc

We consider a linear system of the form

$$\boldsymbol{L}\boldsymbol{x} = \boldsymbol{b}, \tag{23}$$

where $\boldsymbol{L}$ is defined according to the bilinear form of the discretized weak problem, $\boldsymbol{x}$ is the solution vector and $\boldsymbol{b}$ is defined according to the linear form of the discretized weak problem. Multigrid methods try to find the solution to the system in Equation 23 from the solutions to a series of smaller problems. The main assumptions are that high-frequency errors can be effectively eliminated on finer disretizations and smooth errors can be adequately represented on coarser disretizations. In geometric multigrid methods, the problem hierarchy is generated as a set of $n$ progressively coarser discretizations, $\Omega_{e,h}^n \ldots \Omega_{e,h}^1$, where $\Omega_{e,h}^n$ corresponds to the finest problem and $\Omega_{e,h}^1$ is referred to as the base grid and corresponds to the coarsest problem. We discuss the main components of the geometric multigrid method below, which in addition to the grid hierarchy, include transfer operators, smoothers and the base solver.

We make use of space tree data structures [12] for the spatial discretization of the computational domain. Mesh refinement is performed through bisection in the physical space and corresponds to the division of a coarse cell into its $2^d$ children, where $d$ is the dimension. We impose the usual 2:1 balance such that neighboring cells are apart by at most one level of refinement. As a consequence of AMR, hanging nodes are introduced in the discretized system. We handle hanging nodes as constraints and remove them from the global system of equations. We produce the grid hierarchy for geometric multigrid top-down according to Algorithm 1. This process is repeated $n$ times to generate a set of $n$ nested spaces. Transfer operators are responsible for the flow of information between grid levels, i.e., restriction from

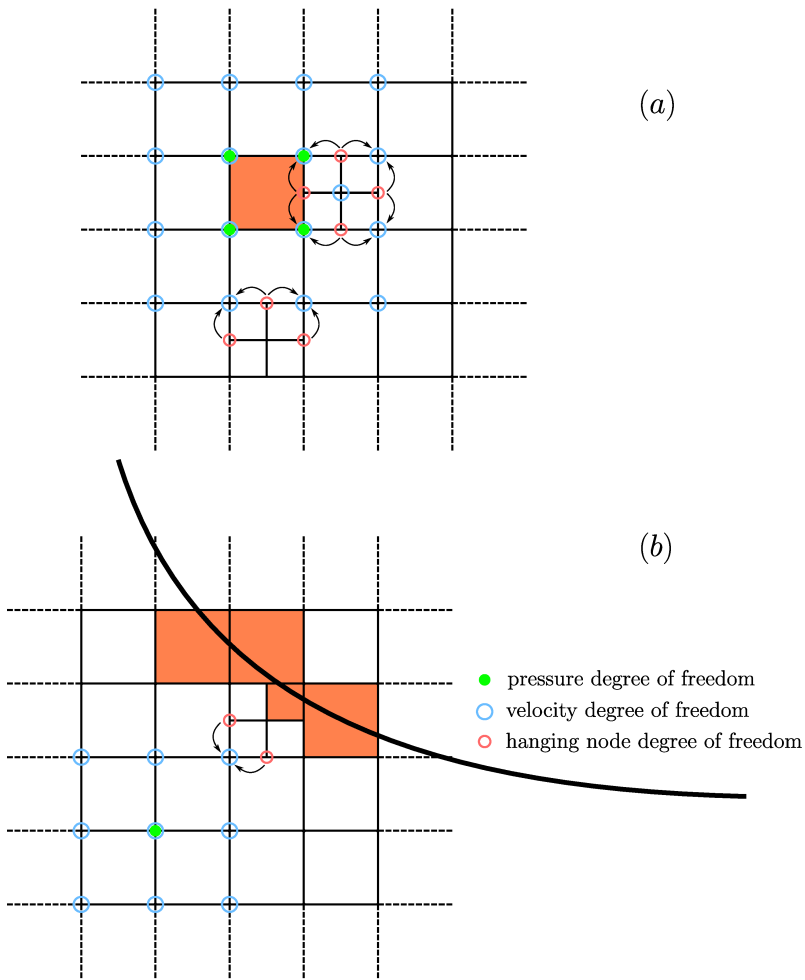

Figure 3: (a) Subdomains in the cell-based smoother. The designated degrees of freedom belong to the subdomain of the shaded cell. A subdomain is generated for every cell in the domain in the same manner. And (b) subdomains in the cutcell-based smoother. A cell-based subdomain is generated for cells cut by the physical boundary (shaded cells in (b)) according to (a), while a node-based subdomain is generated for all pressure nodes that do not appear in any of the cell-based subdomains. A node-based subdomain is illustrated for the designated pressure node in (b). The arrows demonstrate how hanging nodes are implicitly included in the subdomains through their non-hanging counterparts

```
input  : Ω_{e,h}^l
output: Ω_{e,h}^{l-1}

Ω_{e,h}^{l-1} ← Ω_{e,h}^l ;
// r_c is the refinement level of cell c
r_max ← max(r_c ∈ Ω_{e,h}^{l-1}) ;
for c ∈ Ω_{e,h}^{l-1} do
    if r_c == r_max then
        // This replaces c and all its siblings with their parent
        coarsen c ;
    end
end
```
Apply 2:1 balance on $\Omega_{e,h}^{l-1}$ ;

**Algorithm 1:** Generation of the coarse grid $\Omega_{e,h}^{l-1}$ from the fine grid $\Omega_{e,h}^l$

grid level $l$ to $l-1$ and prolongation from grid level $l-1$ to $l$ and can be represented as

$$\begin{aligned} \boldsymbol{w}^{l-1} &= \boldsymbol{R}^l \boldsymbol{w}^l \\ \boldsymbol{w}^l &= \boldsymbol{P}^{l-1} \boldsymbol{w}^{l-1}, \end{aligned} \quad (24)$$

where $\boldsymbol{R}$ and $\boldsymbol{P}$ are the restriction and prolongation operators, respectively, and $\boldsymbol{w}$ is a vector, see also [23, 43].

The success of the multigrid method heavily depends on the smoother. The Vanka smoother [48] has been shown to be an effective choice for the finite difference and finite element discretizations of flow problems. In unfitted finite element methods, it is known that small cut configurations, where the physical boundary intersects a small fraction of a cell lead to severe ill-conditioning of the system matrix. Therefore, the treatment of cut cells in the smoother is an essential aspect of multigrid methods for the finite cell formulation of saddle-point problems. We present two Vanka-type smoothers in the context of Schwarz domain decomposition methods. The smoother operator can be written as

$$\boldsymbol{x}^{k+1} = \boldsymbol{x}^k + \boldsymbol{S}(\boldsymbol{b} - \boldsymbol{L}\boldsymbol{x}^k) \quad (25)$$

where $\boldsymbol{S}$ is the smoother operator.

We formulate two Vanka-type smoothers for the model problems in the context of Schwarz domain decomposition methods. For a grid level $\Omega_{e,h}^l$, let $\mathcal{S}_x := \{x_1, \ldots, x_{n_x}\}$ be the set of all degrees of freedom, and let $\mathcal{S}_u := \{u_1, \ldots, u_{n_u}\} \subset \mathcal{S}_x$ and $\mathcal{S}_p := \{p_1, \ldots, p_{n_p}\} \subset \mathcal{S}_x$ be the sets of velocity and pressure degrees of freedom, respectively, where $\mathcal{S}_x = \mathcal{S}_u \cup \mathcal{S}_p$, $n_x$ is the number of all degrees of freedom, and $n_u$ and $n_p$ are the number of velocity and pressure degrees of freedom, respectively. Define a set of $n_{\text{sd}}$ subdomains $\mathcal{S}_{\text{sd}} := \{\Omega_{\text{sd},1}, \ldots, \Omega_{\text{sd},n_{\text{sd}}}\}$, where subdomain $\Omega_{\text{sd},i}$ spans the space $\mathbb{R}^{n_{\text{sd},i}}$. The general framework of the Schwarz method is to obtain local corrections to the global system based on the solution of the local systems of the subdomains. Given the associated subdomain restriction operators $\boldsymbol{R}_{\text{sd},i}$, the local system of subdomain $\Omega_{\text{sd},i}$ is given by $\boldsymbol{L}_{\text{sd},i} := \boldsymbol{R}_{\text{sd},i} \boldsymbol{L} \boldsymbol{R}_{\text{sd},i}^T$. The restriction operator $\boldsymbol{R}_{\text{sd},i} : \mathbb{R}^{n_x} \to \mathbb{R}^{n_{\text{sd},i}}$ extracts the set of DoFs $\mathcal{S}_{x,\Omega_{\text{sd},i}} := \{x_j \in \mathcal{S}_x \mid x_j \in \Omega_{\text{sd},i}\}$ and its transpose $\boldsymbol{R}_{\text{sd},i}^T : \mathbb{R}^{n_{\text{sd},i}} \to \mathbb{R}^{n_x}$ injects a vector from the local subspace of the subdomain to the global space.

First, we define a cell-based smoother as

$$\boldsymbol{S}_{\text{cell}} = \prod_{i=1}^{n_{\text{sd}}} (\boldsymbol{R}_{\text{sd},i}^{\text{cell}\,T} \boldsymbol{\omega}_i \boldsymbol{L}_{\text{sd},i}^{-1} \boldsymbol{R}_{\text{sd},i}^{\text{cell}}), \quad (26)$$

where $n_{\text{sd}} = n_{\text{c}}$ ($n_{\text{c}}$ is the number of cells in $\Omega_{e,h}$), and $\boldsymbol{\omega}_i$ is in general a diagonal weighting matrix. Let $\mathcal{S}_{p,c_i} := \{p_j \mid p_j \in c_i\}$ be the set of all pressure degrees of freedom that appear in cell $c_i$. Then, subdomain $\Omega_{sd,i}$ contains $\mathcal{S}_{p,c_i}$ and all the velocity degrees of freedom connected to $\mathcal{S}_{p,c_i}$, where connection simply refers to the coupling between degrees of freedom in the discretized bilinear form. It is algebraically equivalent to the degrees of freedom associated with the nonzero entries in $(\boldsymbol{B}^T)_j$ for all $j$ with $p_j \in \mathcal{S}_{p,c_i}$ (see Equation 13). The cell-based subdomain is illustrated in Figure 3(a).

We define a variant of the cell-based smoother next. In this variant, a cell-based subdomain is generated for every cutcell as explained above. The pressure degrees of freedom that do not appear in any of the cell-based subdomains form a set $\mathcal{S}_{p,\text{nb}} := \{p_i \mid p_i \notin \mathcal{S}_{\text{cc}}\}$. For all $p_i \in \mathcal{S}_{p,\text{nb}}$, a node-based subdomain is generated that contains $p_i$ along with all the velocity DoFs connected to it, see Figure 3(b). We denote the smoother operator of this variant as

$$\boldsymbol{S}_{\text{cutcell}} = \prod_{i=1}^{n_{\text{sd}}} (\boldsymbol{R}_{i,\text{sd}}^{\text{cutcell}^T} \boldsymbol{\omega}_i \boldsymbol{L}_i^{-1} \boldsymbol{R}_{i,\text{sd}}^{\text{cutcell}}). \tag{27}$$

where $n_{\text{sd}} = n_{\text{cc}} + |\mathcal{S}_{p,\text{nb}}|$, and $|\mathcal{S}|$ is the cardinality of set $\mathcal{S}$. Both smoothers can be interpreted as multiplicative Schwarz methods, and we note that, conceptionally, the above ideas can be applied to other saddle-point problems.

*Remark* 8. The majority of subdomains in the cutcell-based smoother have a smaller dimension than their counterparts in the cell-based smoother since the physical domain typically intersects only a small fraction of all the cells in the computational domain, i.e., $n_{\text{cc}} \ll n_c$. Therefore, the computational cost of $\boldsymbol{S}_{\text{cutcell}}$ is expected to be noticeably lower than $\boldsymbol{S}_{\text{cell}}$. On the other hand, larger subdomains usually increase the effectiveness of the smoother. The performance and computational cost of both smoothers are analyzed in Section 4.

## 4 Numerical benchmarks

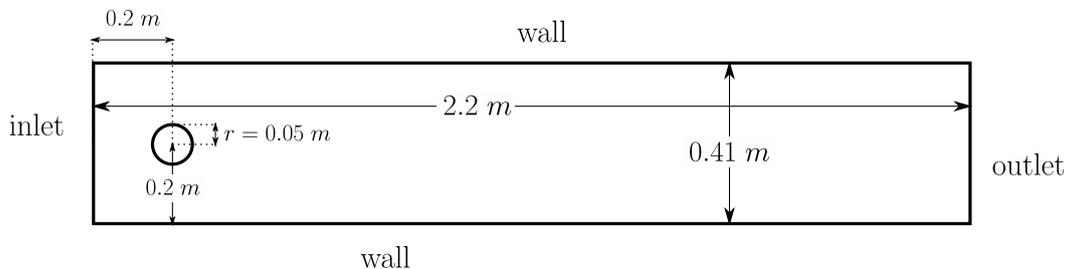

Figure 4: The geometry and boundary conditions of the cylinder flow benchmark. Note that the outlet boundary condition can equivalently be described using a homogeneous Neumann boundary condition. The inlet boundary condition is defined in Equation 28

We consider a numerical benchmark to investigate the performance of the developed geometric multigrid method for benchmark problems. We consider a channel with dimensions $2.2\ m \times 0.41\ m$ with a cylinder with a diameter of $0.1\ m$ near the inflow, see also [44]. The geometry of the channel is shown in Figure 4. The parabolic inflow boundary condition is described as

$$\begin{aligned} u_x(0,y) &= \frac{4\bar{u}y(H-y)}{H^2} & y \in [0,H], \\ u_y(0,y) &= 0 & y \in [0,H], \end{aligned} \tag{28}$$

where $\bar{u} = 0.3\ m/s$. The kinematic viscosity of the fluid is $\eta = \frac{1}{1000}\ m^2/t$.

Geometric multigrid is used both as a solver and as a preconditioner in a Krylov subspace solver. We note that Krylov accelerators are expected to improve the performance of geometric multigrid methods in general; however, in order to better analyze the performance of the developed multigrid method for the systems of equations that arise from the model problems, we use multigrid also as a standalone fixed-point iterator in order to exclude the added effect of Krylov accelerators. A V-cycle is used in all examples. We consider the two smoother variants introduced in Section 3 in order to provide insight into the performance of the geometric multigrid method with respect to the choice of smoothers and the required computational cost. The smoothers differ solely in the Schwarz subdomains. A subdomain is generated for each cell in the first variant. On the other hand, cell-based subdomains are only generated for pressure nodes that appear in a cutcell and node-based subdomains are employed for pressure nodes that do not appear in any cutcell in the second variant. All examples are executed on an Intel Xeon Skylake Gold 6148 CPU with 192 GB of main memory using an in-house C++ implementation with `p4est` [12] and `PETSc` [3] for mesh manipulation and some linear algebra components, respectively, and `Paraview` [1] for visualization postprocessing.

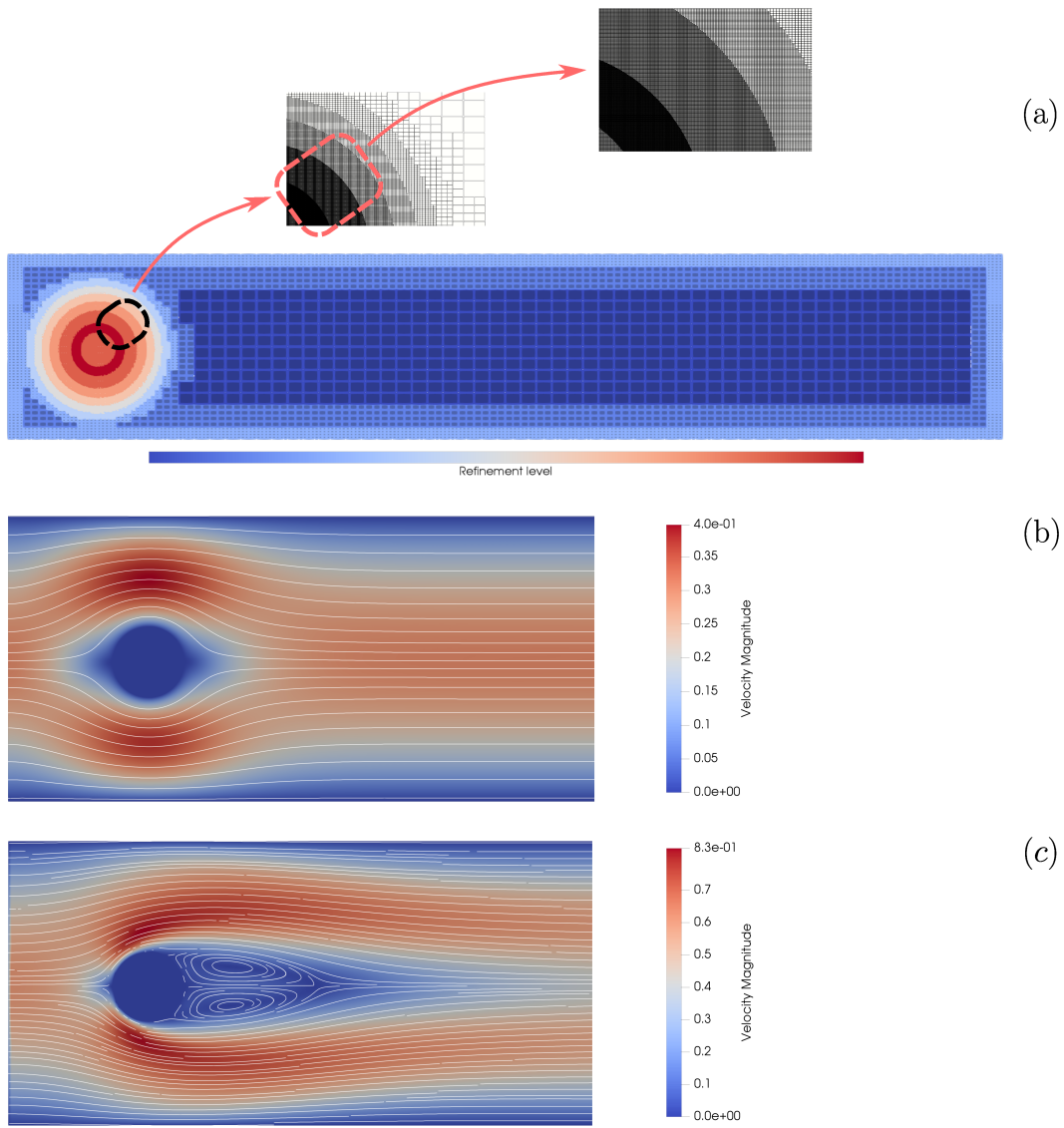

Figure 5: The cylinder flow benchmark: (a) adaptive refinement of the computational mesh towards the boundaries and the cylinder, where the red and blue spectra indicate lower and higher refinement levels, respectively, (b) the velocity magnitude and streamlines for the Stokes problem and (c) the velocity magnitude and streamlines for the Navier-Stokes problem. Note that only the left part of the domain is shown in (b) and (c)

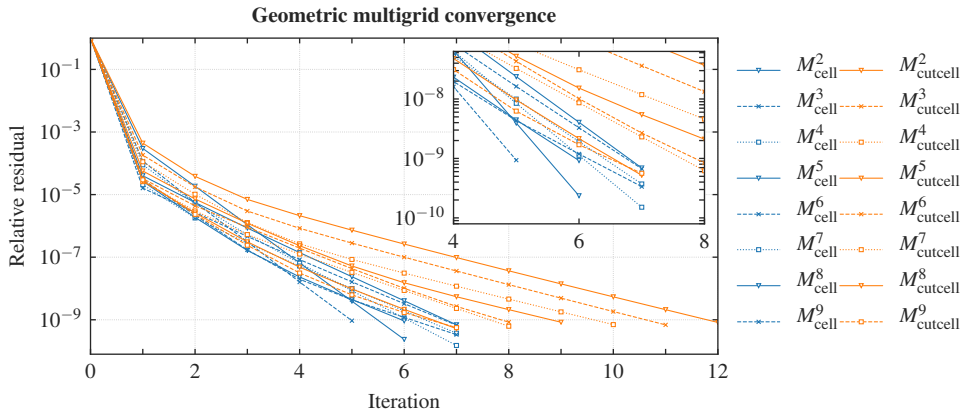

Figure 6: The convergence of the geometric multigrid solver for the Stokes flow in the cylinder flow benchmark problem. Geometric multigrid is used as a solver. The grid hierarchy given in Table 1 is used for the mesh study

Table 1: The grid hierarchy used in the convergence study of the cylinder flow benchmark. $n_c$ and $n_{\mathrm{DoF}}$ are the number of cells and the number of degrees of freedom, respectively

| Grid | $n_c$ | $n_{\mathrm{DoF}}$ |
| --- | ---: | ---: |
| $\Omega_h^9$ | 881,155 | 2,633,358 |
| $\Omega_h^8$ | 525,421 | 1,571,277 |
| $\Omega_h^7$ | 194,485 | 580,737 |
| $\Omega_h^6$ | 67,090 | 199,962 |
| $\Omega_h^5$ | 23,857 | 71,082 |
| $\Omega_h^4$ | 10,054 | 30,138 |
| $\Omega_h^3$ | 5,767 | 17,535 |
| $\Omega_h^2$ | 2,581 | 8,040 |
| $\Omega_h^1$ | 1,024 | 3,315 |

We study the Stokes problem and the Navier-Stokes problem separately on the cylinder flow benchmark. In each case, the convergence of the solver is analyzed through a mesh study using a progressively finer set of problems. A combination of uniform and adaptive

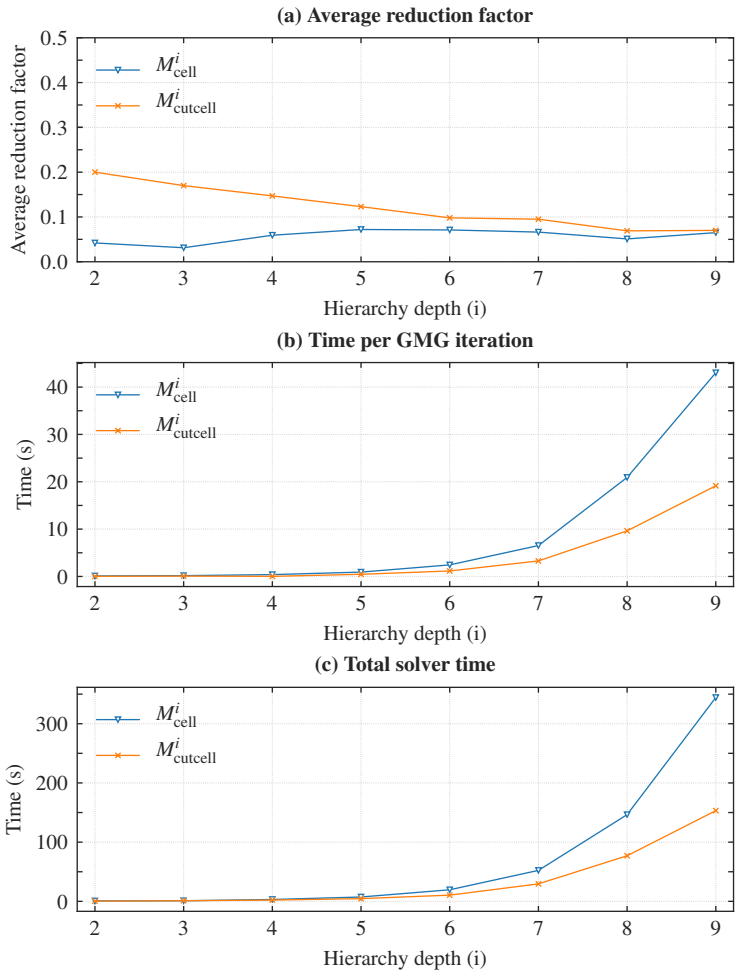

Figure 7: The performance of the two smoother variants within the geometric multigrid solver for the Stokes flow in the cylinder flow benchmark problem: (a) the average reduction rate of the solver, (b) time per geometric multigrid iteration and (c) total solver runtime. The grid hierarchy given in Table 1 is used

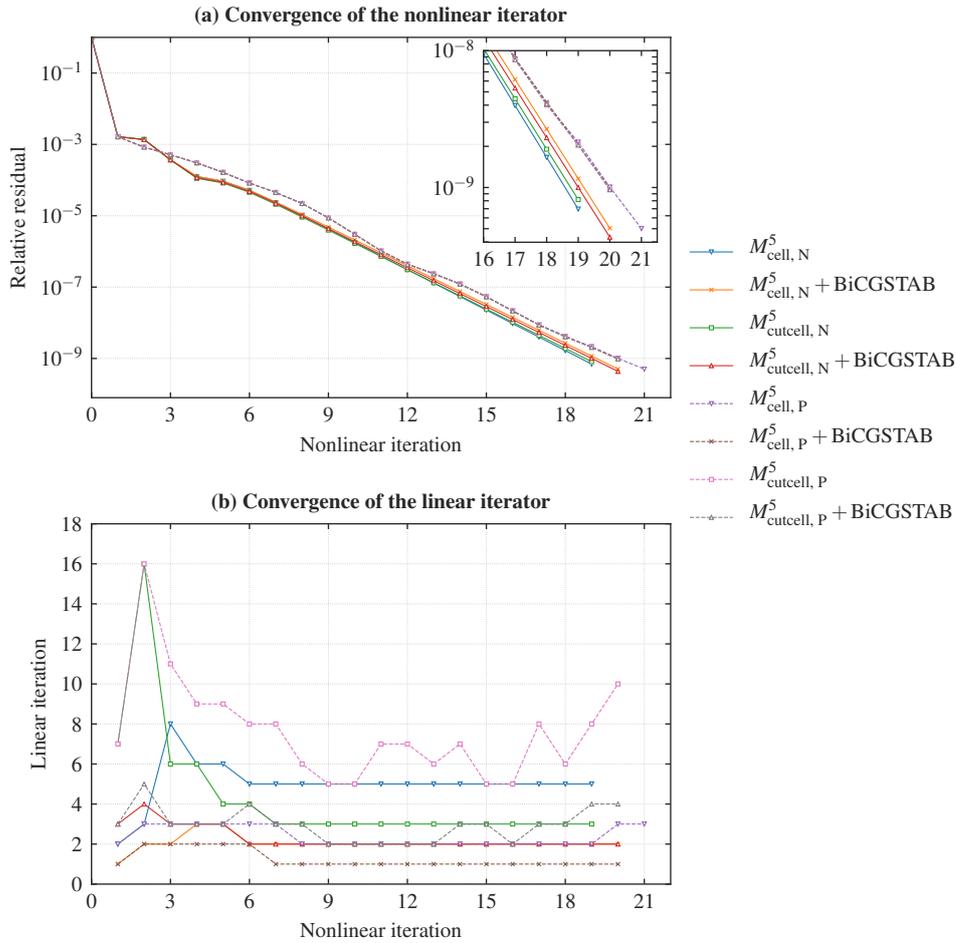

Figure 8: The convergence of the nonlinear iterator using Newton and Picard linearization (a) and the convergence of the linear iterator at each linearization point using geometric multigrid as a solver and as a preconditioner in the BiCGSTAB method (b). The Newton and Picard linearization methods are designated by $N$ and $P$, respectively. In the absence of explicit indication, the geometric multigrid is used as a solver

refinement steps is employed to produce the computational meshes. Adaptive refinement is performed near the boundaries of the channel and around the cylinder as shown in Figure 5(a). The generated grid hierarchy is given in Table 1. $\Omega_h^1$ is used as the coarse grid for all problems, i.e., finer problems employ a larger number of grid levels. $M_{\boldsymbol{S}}^i$ designates a problem with $\Omega_h^i$ as the fine grid and consequently $i$ levels of grid hierarchy. $\boldsymbol{S}$ denotes the smoother operator and is either cell for the cell-based smoother or cutcell for the cutcell-based smoother. In addition to adaptive refinement, eight steps of adaptive integration are used for cutcells.

First, we consider the Stokes problem. Geometric multigrid with three steps of pre- and post-smoothing is used as a solver. A damping factor of $2/3$ is used for both smoothers. A tolerance of $10^{-9}$ on the relative residual is used as the convergence criterion. The velocity profile of the fluid around the cylinder is shown in Figure 5(b). The convergence of the multigrid solver is shown in Figure 6. The multigrid method demonstrates rapid convergence for all problems $M^i$, $i = 2, \ldots, 9$. It can be observed that roughly the same number of iterations are required for the convergence of problems with different fine grids and consequently different depths of grid hierarchy, which illustrates the independence of the solver from the problem size and demonstrates its capability to employ a deep hierarchy.

The performance of the cell-based and cutcell-based smoothers, including average reduction factor, time per GMG iteration and total solver time is shown in Figure 7. The cell-based smoother consistently achieves a lower reduction factor. Consequently, the required number of iterations for convergence by the cell-based smoother is typically smaller as seen in Figure 6. However, the difference between the average reduction factors of the smoothers is not substantial, especially on finer grids with deeper grid hierarchies (see Figure 7(a)). On the other hand, the cell-based smoother is more computationally intensive, which is reflected by its per-iteration runtime being noticeably larger than its cutcell-based counterpart as shown in Figure 7(b). As a result, the total solver runtime is consistently smaller for the cutcell-based smoother, despite its ostensible disadvantage with regards to shear residual reduction rate (see Figure 7(c)).

We turn our focus to the Navier-Stokes problem next. The range of solution strategies, including the linearization method, the choice of the Krylov accelerator, the geometric multigrid settings, etc. gives rise to a myriad of possible configurations. We try to discuss the most important aspects of such configurations in the following. We consider geometric multigrid both as a standalone solver and as a preconditioner the BiCGSTAB method [41] for the solution of the linearized system. At each linearization point, the linear iteration is stopped when the relative residual is reduced by a factor of $10^2$. The initial guess is the zero vector. The velocity profile of the fluid near the cylinder is shown in Figure 5(c).

We first examine the presented linearization methods, namely Newton and Picard linearization (see Section 2). Three pre- and post-smoothing steps are used in the geometric multigrid method. A damping factor of $2/3$ for the cell-based smoother and harmonic weighting for the cutcell-based smoother are used. The convergence of the nonlinear iteration as well as the required number of linear iterations at each linearization point is shown in Figure 8 for problem $M^5$. Although the Picard iteration leads to a slight increase in the number of nonlinear iterations, the overall convergence of the Newton and Picard iterations are similar for the presented benchmark. On the other hand, the BiCGSTAB method reduces the required number of linear iterations in all configurations (see Figure 8(b)) and, in the majority of cases, only one or two iterations are sufficient to achieve the required relative reduction of the residual.

*Remark* 9. We note that the reduction of the relative residual by a larger factor, i.e., more linear iterations at the linearization points does not significantly affect the overall convergence of the nonlinear iteration for the presented benchmark. Therefore, in favor of minimizing computational cost, we reduce the relative residual by a factor of $10^2$.

*Remark* 10. Despite the apparent similarity between the convergence of Newton and Picard iterations for the presented benchmark, it is expected that the Newton iteration achieve a more rapid rate of convergence in general, given that the iteration is sufficiently close to the solution. In this regard, we also note that the initial guess for the solution of the Navier-Stokes system becomes increasingly important for higher Reynold numbers. In such cases, one standard strategy is to solve a series of problems with increasing Reynolds numbers and use the solution to the Stokes system as the initial guess for the first problem.

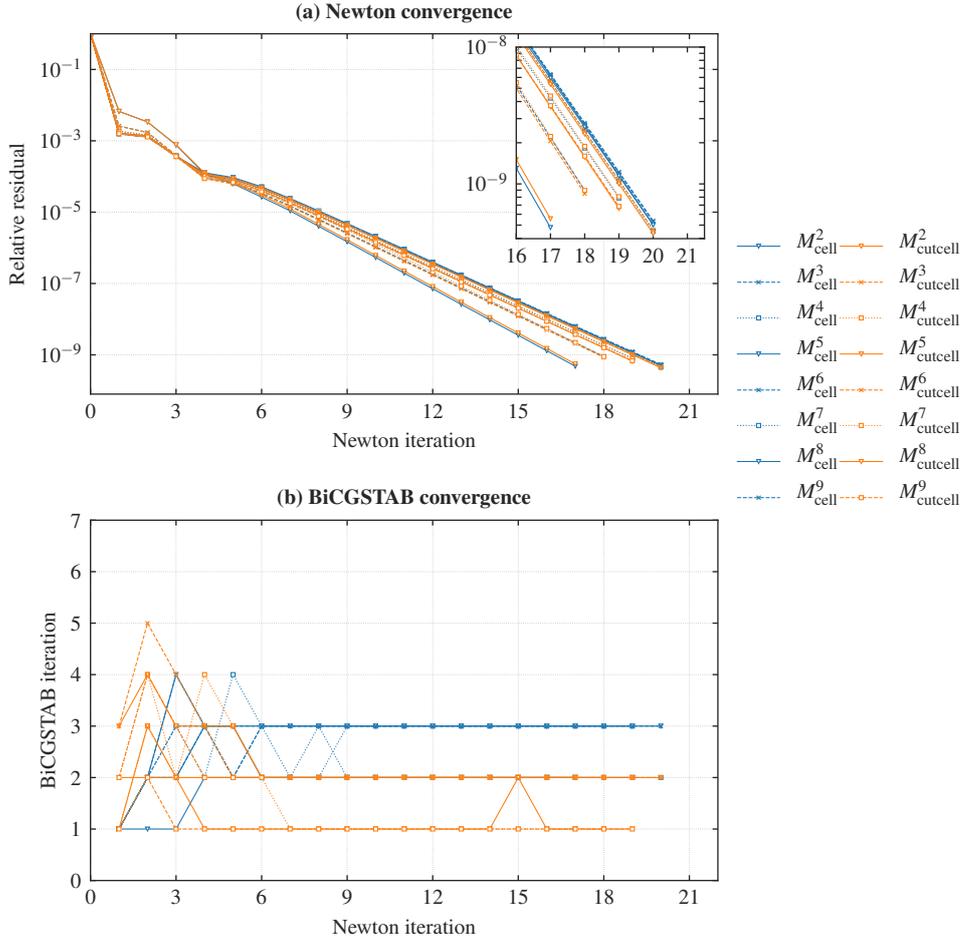

Figure 9: The convergence of the Newton solver for the Navier-Stokes problem in the cylinder flow benchmark (a). A BiCGSTAB solver with geometric multigrid as its preconditioner is used to solve the linearized system at each step. The required number of BiCGSTAB iterations to achieve a relative convergence of $r_{\text{rel}} = 10^{-2}$ at each Newton linearization point is shown in (b)

We analyze the convergence of the geometric multigrid method for the Navier-Stokes problem through a mesh study next. We employ the grid hierarchy in Table 1. The system is linearized using the Newton's method, and the geometric multigrid method is used as a preconditioner in the BiCGSTAB method to solve the linearized systems. The convergence of the nonlinear iterator as well as the number of required BiCGSTAB iterations at each linearization point is shown in Figure 9. It can be seen that the variations in the number of required nonlinear iterations for different problem sizes are small. As expected, the overall convergence of the Newton iteration is not meaningfully affected by the choice of smoother in the geometric multigrid method. We note that minor differences can occur since the linearized systems are not solved to the same exact residual—the relative residual is reduced by a factor of $10^2$. Furthermore, the number of required BiCGSTAB iterations are roughly constant and have an inclination to reach a plateau—one to three iterations in the presented mesh study—after the first few nonlinear iterations (see Figure 9(b)). Three steps of pre- and post-smoothing are used for all problems except $M_{cutcell}^i$ $i = 7, 8, 9$, where six to nine smoothing steps were required. We note that three smoothing steps would be sufficient if a slightly finer mesh were used as the course grid in these problems. It is therefore induced that the cell-based smoother is more effective per iteration for very deep grid hierarchies for the Naiver-Stokes problem. Nevertheless, the solver demonstrates excellent convergence independently of the problem size and hierarchy depth for the presented benchmark.

# 5 Conclusions

Unfitted finite element methods circumvent the generation of a boundary-conforming mesh, which is a bottleneck in classical finite element methods. On the other hand, the discrete systems arising from such methods, including the finite cell method, are typically ill-conditioned that in turn renders the efficient solution of such problems a challenging task. We employed an adaptive geometric multigrid method for the solution of the mixed finite cell formulation of saddle-point problems in the context of the Stokes and Navier-Stokes equations. We presented two smoothers tailored to the finite cell formulation of saddle-point problems and studied the convergence of the multigrid method through a numerical benchmark. Results indicate that the multigrid method both as a solver and as a preconditioner in a Krylov subspace solver can effectively solve the linear(ized) systems independently of the problem size. It is also seen that the multigrid method is robust with respect to the grid hierarchy and deep grid hierarchies allow the solution to large problems to be obtained using comparatively small problems. The cell-based smoother is more effective per iteration as compared to the cutcell-based smoother. On the other hand, the cutcell-based smoother is more favorable in terms of computational cost compared to its cell-based counterpart and achieves a faster total solution time in the presented benchmark.

**Acknowledgments** Financial support was provided by the German Research Foundation (*Deutsche Forschungsgemeinschaft, DFG*) in the framework of subproject C4 of the Collaborative Research Center SFB 837 *Interaction Modeling in Mechanized Tunneling*. This support is gratefully acknowledged.

# References

[1] Ahrens, J., Geveci, B., Law, C.: ParaView: An End-User Tool for Large Data Visualization, Visualization Handbook. Elsevier (2005)

[2] Babuška, I.: The finite element method with penalty. Mathematics of computation **27**(122), 221–228 (1973)

[3] Balay, S., Gropp, W.D., McInnes, L.C., Smith, B.F.: Efficient management of parallelism in object oriented numerical software libraries. In: E. Arge, A.M. Bruaset, H.P. Langtangen (eds.) Modern Software Tools in Scientific Computing, pp. 163–202. Birkhäuser Press (1997)

[4] Belytschko, T., Moës, N., Usui, S., Parimi, C.: Arbitrary discontinuities in finite elements. International Journal for Numerical Methods in Engineering **50**(4), 993–1013 (2001)

[5] Braess, D., Sarazin, R.: An efficient smoother for the stokes problem. Applied Numerical Mathematics **23**(1), 3–19 (1997)

[6] Bramble, J.H., Pasciak, J.E.: A preconditioning technique for indefinite systems resulting from mixed approximations of elliptic problems. Mathematics of Computation **50**(181), 1–17 (1988)

[7] Bramble, J.H., Pasciak, J.E., Vassilev, A.T.: Analysis of the inexact uzawa algorithm for saddle point problems. SIAM Journal on Numerical Analysis **34**(3), 1072–1092 (1997)

[8] Burman, E.: Ghost penalty. Comptes Rendus Mathematique **348**(21), 1217 – 1220 (2010). DOI https://doi.org/10.1016/j.crma.2010.10.006. URL http://www.sciencedirect.com/science/article/pii/S1631073X10002827

[9] Burman, E., Claus, S., Hansbo, P., Larson, M.G., Massing, A.: Cutfem: discretizing geometry and partial differential equations. International Journal for Numerical Methods in Engineering **104**(7), 472–501 (2015)

[10] Burman, E., Hansbo, P.: Fictitious domain finite element methods using cut elements: I. a stabilized lagrange multiplier method. Computer Methods in Applied Mechanics and Engineering **199**(41-44), 2680–2686 (2010)

[11] Burman, E., Hansbo, P.: Fictitious domain finite element methods using cut elements: Ii. a stabilized nitsche method. Applied Numerical Mathematics **62**(4), 328–341 (2012)

[12] Burstedde, C., Wilcox, L.C., Ghattas, O.: p4est: Scalable algorithms for parallel adaptive mesh refinement on forests of octrees. SIAM Journal on Scientific Computing **33**(3), 1103–1133 (2011)

[13] Dolbow, J., Harari, I.: An efficient finite element method for embedded interface problems. International journal for numerical methods in engineering **78**(2), 229–252 (2009)

[14] Düster, A., Parvizian, J., Yang, Z., Rank, E.: The finite cell method for three-dimensional problems of solid mechanics. Computer methods in applied mechanics and engineering **197**(45-48), 3768–3782 (2008)

[15] Elman, H., Howle, V.E., Shadid, J., Shuttleworth, R., Tuminaro, R.: Block preconditioners based on approximate commutators. SIAM Journal on Scientific Computing **27**(5), 1651–1668 (2006)

[16] Elman, H.C.: Multigrid and krylov subspace methods for the discrete stokes equations. International journal for numerical methods in fluids **22**(8), 755–770 (1996)

[17] Elman, H.C., Golub, G.H.: Inexact and preconditioned uzawa algorithms for saddle point problems. SIAM Journal on Numerical Analysis **31**(6), 1645–1661 (1994)

[18] Embar, A., Dolbow, J., Harari, I.: Imposing dirichlet boundary conditions with nitsche's method and spline-based finite elements. International journal for numerical methods in engineering **83**(7), 877–898 (2010)

[19] Fernández-Méndez, S., Huerta, A.: Imposing essential boundary conditions in meshfree methods. Computer methods in applied mechanics and engineering **193**(12-14), 1257–1275 (2004)

[20] Flemisch, B., Wohlmuth, B.I.: Stable lagrange multipliers for quadrilateral meshes of curved interfaces in 3d. Computer Methods in Applied Mechanics and Engineering **196**(8), 1589–1602 (2007)

[21] Gjesdal, T., Lossius, M.E.H.: Comparison of pressure correction smoothers for multigrid solution of incompressible flow. International Journal for Numerical Methods in Fluids **25**(4), 393–405 (1997)

[22] Glowinski, R., Kuznetsov, Y.: Distributed lagrange multipliers based on fictitious domain method for second order elliptic problems. Computer Methods in Applied Mechanics and Engineering **196**(8), 1498–1506 (2007)

[23] Hackbusch, W.: Multi-grid methods and applications. Springer Series in Computational Mathematics **4** (1985)

[24] Hansbo, A., Hansbo, P.: An unfitted finite element method, based on nitsches method, for elliptic interface problems. Computer methods in applied mechanics and engineering **191**(47-48), 5537–5552 (2002)

[25] Hoang, T., Verhoosel, C.V., Auricchio, F., van Brummelen, E.H., Reali, A.: Mixed isogeometric finite cell methods for the stokes problem. Computer Methods in Applied Mechanics and Engineering **316**, 400–423 (2017)

[26] Jomo, J.N., de Prenter, F., Elhaddad, M., D'Angella, D., Verhoosel, C.V., Kollmannsberger, S., Kirschke, J.S., Nübel, V., van Brummelen, E., Rank, E.: Robust and parallel scalable iterative solutions for large-scale finite cell analyses. Finite Elements in Analysis and Design **163**, 14–30 (2019)

[27] Jomo, J.N., Zander, N., Elhaddad, M., Özcan, A., Kollmannsberger, S., Mundani, R.P., Rank, E.: Parallelization of the multi-level hp-adaptive finite cell method. Computers & Mathematics with Applications **74**(1), 126–142 (2017)

[28] Kay, D., Loghin, D., Wathen, A.: A preconditioner for the steady-state navier–stokes equations. SIAM Journal on Scientific Computing **24**(1), 237–256 (2002)

[29] Kelley, C.T.: Iterative methods for linear and nonlinear equations. SIAM (1995)

[30] Klawonn, A., Pavarino, L.F.: Overlapping schwarz methods for mixed linear elasticity and stokes problems. Computer Methods in Applied Mechanics and Engineering **165**(1-4), 233–245 (1998)

[31] Klawonn, A., Pavarino, L.F.: A comparison of overlapping schwarz methods and block preconditioners for saddle point problems. Numerical linear algebra with applications **7**(1), 1–25 (2000)

[32] Larin, M., Reusken, A.: A comparative study of efficient iterative solvers for generalized stokes equations. Numerical linear algebra with applications **15**(1), 13–34 (2008)

[33] Nitsche, J.: Über ein variationsprinzip zur lösung von dirichlet-problemen bei verwendung von teilräumen, die keinen randbedingungen unterworfen sind. In: Abhandlungen aus dem mathematischen Seminar der Universität Hamburg, vol. 36, pp. 9–15. Springer (1971)

[34] Parvizian, J., Düster, A., Rank, E.: Finite cell method. Computational Mechanics **41**(1), 121–133 (2007)

[35] Patankar, S.V., Spalding, D.B.: A calculation procedure for heat, mass and momentum transfer in three-dimensional parabolic flows. In: Numerical prediction of flow, heat transfer, turbulence and combustion, pp. 54–73. Elsevier (1972)

[36] Pavarino, L.F.: Indefinite overlapping schwarz methods for time-dependent stokes problems. Computer methods in applied mechanics and engineering **187**(1-2), 35–51 (2000)

[37] Peters, J., Reichelt, V., Reusken, A.: Fast iterative solvers for discrete stokes equations. SIAM journal on scientific computing **27**(2), 646–666 (2005)

[38] de Prenter, F., Verhoosel, C., van Brummelen, E.: Preconditioning immersed isogeometric finite element methods with application to flow problems. Computer Methods in Applied Mechanics and Engineering **348**, 604–631 (2019)

[39] de Prenter, F., Verhoosel, C.V., van Brummelen, E., Evans, J., Messe, C., Benzaken, J., Maute, K.: Multigrid solvers for immersed finite element methods and immersed isogeometric analysis. Computational Mechanics pp. 1–32 (2019)

[40] Rusten, T., Winther, R.: A preconditioned iterative method for saddlepoint problems. SIAM Journal on Matrix Analysis and Applications **13**(3), 887–904 (1992)

[41] Saad, Y.: Iterative methods for sparse linear systems, vol. 82. siam (2003)

[42] Saberi, S., Vogel, A., Meschke, G.: Parallel finite cell method with adaptive geometric multigrid. In: European Conference on Parallel Processing, pp. 578–593. Springer (2020)

[43] Sampath, R.S., Biros, G.: A parallel geometric multigrid method for finite elements on octree meshes. SIAM Journal on Scientific Computing **32**(3), 1361–1392 (2010)

[44] Schäfer, M., Turek, S., Durst, F., Krause, E., Rannacher, R.: Benchmark computations of laminar flow around a cylinder. In: Flow simulation with high-performance computers II, pp. 547–566. Springer (1996)

[45] Schillinger, D., Ruess, M.: The finite cell method: A review in the context of higher-order structural analysis of cad and image-based geometric models. Archives of Computational Methods in Engineering **22**(3), 391–455 (2015)

[46] Silvester, D., Elman, H., Kay, D., Wathen, A.: Efficient preconditioning of the linearized navier–stokes equations for incompressible flow. Journal of Computational and Applied Mathematics **128**(1-2), 261–279 (2001)

[47] Silvester, D., Wathen, A.: Fast iterative solution of stabilised stokes systems part ii: using general block preconditioners. SIAM Journal on Numerical Analysis **31**(5), 1352–1367 (1994)

[48] Vanka, S.P.: Block-implicit multigrid solution of navier-stokes equations in primitive variables. Journal of Computational Physics **65**(1), 138–158 (1986)

[49] Verfürth, R.: A multilevel algorithm for mixed problems. SIAM journal on numerical analysis **21**(2), 264–271 (1984)

[50] Wittum, G.: Multi-grid methods for stokes and navier-stokes equations. Numerische Mathematik **54**(5), 543–563 (1989)

[51] Zhu, T., Atluri, S.: A modified collocation method and a penalty formulation for enforcing the essential boundary conditions in the element free galerkin method. Computational Mechanics **21**(3), 211–222 (1998)